\documentstyle[hyperref]{article}

\newfont{\cirilrm}{wncyr10 scaled 1000}
\newfont{\cirilbf}{wncyb10 scaled 1000}
\newfont{\cirilsf}{wncyss10 scaled 1000}
\newfont{\cirilit}{wncyi10 scaled 1000}
\newfont{\cirilsc}{wncysc10 scaled 1000}

\newfont{\newit}{cmfi10 scaled 1200}
\newfont{\newsit}{cmfi10 scaled 1000}
\newfont{\newsmit}{cmfi10 scaled 800}
\newfont{\newssmit}{cmfi10 scaled 600}

\newfont{\newmat}{msbm10 scaled 1000}
\newfont{\newmatmb}{msbm10 scaled 800}
\newfont{\newmatsm}{msbm10 scaled 600}

\def\ppsm#1{\mbox{\newsmit #1}}
\def\ppssm#1{\mbox{\newssmit #1}}

\newcommand{\z}{\symbol{'31}}

\newcommand{\sh}{\symbol{'170}}

\newcommand{\ch}{\symbol{'161}}

\newcommand{\ja}{\symbol{'37}}

\newcommand{\Ju}{\symbol{'20}}

\newcommand{\ijj}{\symbol{'32}}

\newcommand{\ii}{\symbol{'171}}

\newcounter{supersection}[section]
\newtheorem{th}[supersection]{Theorem}
\newtheorem{lm}[supersection]{Lemma}
\newtheorem{no}[supersection]{Note}
\newtheorem{re}[supersection]{Remark}
\newtheorem{co}[supersection]{Corollary}

\def\bibname{\textbf{REFERENCES}}
\def\thebibliography#1{\paragraph*{\uppercase{\bibname}}\list
{[\arabic{enumi}]}{\settowidth\labelwidth{[#1]}\leftmargin\labelwidth
\advance\leftmargin\labelsep\usecounter{enumi}}
\def\newblock{\hskip .11em plus .33em minus .07em}
\sloppy\clubpenalty4000\widowpenalty4000
\sfcode`\.=1000\relax}

\def\Dj{D{\hspace{-.75em}\raisebox{.3ex}{-}\hspace{.4em}}}

\def\stop{\mbox{\footnotesize {\vrule width 6pt height 6pt}}}

\begin{document}

\thispagestyle{plain}


\noindent $\;$ %

\vspace*{10.00 mm}

\centerline{\large \rm SOME CONSIDERATIONS IN CONNECTION}

\smallskip

\centerline{\large \rm WITH ALTERNATING KUREPA'S FUNCTION}

\medskip

\bigskip
\centerline{\it Branko J. Male\v sevi\' c${\,}^{1}$}
\footnotetext{$\!\!\!\!\!\!\!\!{}^{1}$Research partially supported by the MNTRS, Serbia, Grant No. 144020.}

\begin{center}
${}^{1}$Faculty of Electrical Engineering \\
University of Belgrade, 11000 Belgrade, Serbia \\
E-mail: malesh@eunet.yu
\end{center}

\noindent
{\small {\bf ABSTRACT.}
In this paper we consider the functional equation for alternating factorial sum and some of its particular solutions
{\normalsize (}alternating Kurepa's function $A(z)$ from \cite{Petojevic_02} and function $A_{1}(z)${\normalsize )}.
We determine an extension of domain of functions $A(z)$ and $A_{1}(z)$ in the sense of the principal value at point
\cite{Slavic_73}, \cite{MijajlovicMalesevic_07}. Using the methods from \cite{Slavic_73} and \cite{Malesevic_03}
we give a new representation of alternating Kurepa's function $A(z)$, which is an analog of Slavi\' c's representation
of Kurepa's function $K(z)$ \cite{Slavic_73}, \cite{Marichev_83}. Also, we consider some representations of functions
$A(z)$ and $A_{1}(z)$ via incomplete gamma function and we consider differential transcendency of previous functions~too.}

\medskip
\leftline{KEYWORDS: Gamma function, Kurepa's function, Casimir energy.}

\medskip
\leftline{MSC (2000): 30E20, 11J91, 81V99.}

\bigskip

\bigskip
\noindent
\section{\bf \boldmath \hspace*{-7.0 mm}
1 The functional equation for alternating factorial sum and its particular solutions}

The main object of consideration in this paper is the functional equation for
alternating factorial sum
\begin{equation}
\label{A_FE_1}
A(z) + A(z-1) = \Gamma(z+1),
\end{equation}
with respect to the function $A : \mbox{\newmat{D}} \longrightarrow \mbox{\newmat{C}}$ with domain
$\mbox{\newmat{D}} \subseteq \mbox{\newmat{C}} \backslash \mbox{\newmat{Z}}^{-}$,  where $\Gamma$
is the gamma function, $\mbox{\newmat{C}}$ is the set of complex numbers and $\mbox{\newmat{Z}}^{-}$
is set of negative integer numbers. A solution of functional equation (\ref{A_FE_1})
over the set of natural numbers ($\mbox{\newmat{D}} = \mbox{\newmat{N}}$)
is the function of alternating left factorial $A_n$. R. Guy introduced this function,
in the book \cite{Guy_94} (p.$\:$100), as an alternating sum of factorials
\mbox{$A_n = n! - (n-1)! + \ldots + (-1)^{n-1}1! \, .$}
Let us use the notation
\begin{equation}
\label{A_SUM_1}
A(n) = \displaystyle\sum\limits_{i=1}^{n}{(-1)^{n-i}i!} \; .
\end{equation}

\break

\noindent
Sum (\ref{A_SUM_1}) corresponds to the sequence $A005165$ in \cite{Sloane_06}.
We call the functional equation (\ref{A_FE_1}) {\em the functional equation for
alternating factorial sum}. In consideration which follows we consider two
particular solutions of the functional equation (\ref{A_FE_1}).

\bigskip
\noindent
{\bf \boldmath 1.1$\;$The function $A(z)$} An analytical extension of the function
(\ref{A_SUM_1}) over the set of complex numbers is determined by integral~\cite{Petojevic_02}:
\begin{equation}
\label{A_INT_1}
A(z)
=
\displaystyle\int\limits_{0}^{\infty}{
e^{-t} \displaystyle\frac{t^{z+1}-(-1)^{z}t}{t+1} \: dt},
\end{equation}
which converges for $\mbox{Re} \: z > 0$. For the function $A(z)$ we use the
term {\em alternating Kurepa's function} and it is a solution of the functional
equation (\ref{A_FE_1}). Let us observe that since \mbox{$A(z-1) = \Gamma(z+1) -A(z)$},
it is possible to make analytical continuation of alternating Kurepa's function $A(z)$
for \mbox{$\mbox{Re} \, z \leq 0$}. In that way, the alternating Kurepa's function
$A(z)$ is a meromorphic function with simple poles at $z \!=\! -n$ $(n \!\geq\! 2)$.
At a point $z = \infty$ alternating Kurepa's function has an essential singularity.
Alternating Kurepa's function has the following residues
\begin{equation}
\mathop{\mbox{\rm res}}\limits_{z = -n}{A(z)} =
(-1)^n \displaystyle\sum\limits_{k=0}^{n-2}{\displaystyle\frac{1}{k!}}
\quad (n\!\geq\!2).
\end{equation}
Previous results for alternating Kurepa's function are given according to \cite{Petojevic_02}.

\bigskip
\noindent
{\bf \boldmath 1.2$\;\,$The function $A_{1}(z)\,$}
The functional equation (\ref{A_FE_1}), besides alternating Kurepa's function $A(z)$,
has another solution which is given by the following statement.
\begin{th}
\label{A_lema_4}
Let $\mbox{\newmat{D}} = \mbox{\newmat{C}} \backslash \mbox{\newmat{Z}}$. Then, series
\begin{equation}
\label{Def_A1}
A_{1}(z) = \displaystyle\sum\limits_{n=0}^{\infty}{(-1)^n\Gamma(z+1-n)}
\end{equation}
absolutely converges and it is a solution of the functional equation
{\rm (\ref{A_FE_1})}~over~$\mbox{\newmat{D}}$.
\end{th}

\noindent
{\bf Proof.} Statement of the Theorem is a consequence of the Theorem 1.1
from \cite{Malesevic_03}.~\stop
\begin{re}
Function $A_{1}(z)$, defined by {\rm (\ref{Def_A1})} over $\mbox{\newmat{C}}$,
has poles at integer points $z\!=\!m \in \mbox{\newmat{Z}}$.
\end{re}

\section{\bf \boldmath \hspace*{-7.0 mm}
2 Extending the domain of functions $A(z)$ and $A_{1}(z)$ in the sense of the principal value at point}

Let us observe a possibility of extending the domain of the functions $A(z)$ and $A_{1}(z)$,
in the sense of the principal value at point, over the set of complex numbers. Namely, for
a meromorphic function $f(z)$, on the basis of Cauchy's integral formula, we define
{\em the principal value at point $a$} as follows \cite{Slavic_70}, \cite{Slavic_73}:
\begin{equation}
\label{GAMMA_PV_a}
\mathop{\mbox{\rm p.v.}}\limits_{z = a}{f(z)}
=
\lim\limits_{\rho \rightarrow 0_{+}}{\displaystyle\frac{1}{2 \pi i}
\!\!\!\!\displaystyle\oint\limits_{|z-a|=\rho}{\!\!\!\!\displaystyle\frac{f(z)}{z-a}\,dz}}.
\end{equation}
It is obvious that the principal value at pole $z = a$ exists as a finite complex number
$\mathop{\mbox{res}}\limits_{z=a}{{\Big (}\displaystyle\frac{f(z)}{z-a}}{\Big )}$.
For two meromorphic functions $f_1(z)$ and $f_2(z)$ additivity is true~\cite{Slavic_70}:
\begin{equation}
\mathop{\mbox{\rm p.v.}}\limits_{z = a}{{\Big (}f_{1}(z)+f_{2}(z){\Big )}}
=
\mathop{\mbox{\rm p.v.}}\limits_{z = a}{f_{1}(z)}
+
\mathop{\mbox{\rm p.v.}}\limits_{z = a}{f_{2}(z)}.
\end{equation}
In the paper \cite{Slavic_70} it is proved that multiplicativity
of the principal value does not hold. The following statement is proved in \cite{MijajlovicMalesevic_07}.
\begin{th}
Let $f_{1}(z)$ be a holomorphic function at the point $a$ and let $f_{2}(z)$ be a meromorphic function with
pole of the $m$-th order at the same point $a$.$\;$Then
\begin{equation}
\mathop{\mbox{\rm p.v.}}\limits_{z = a}{ {\Big (} f_{1}(z) \cdot
f_{2}(z) {\Big )} } =
\displaystyle\sum\limits_{k=0}^{m}{\displaystyle\frac{f_{1}^{(k)}(a)}{k!}
\mathop{\mbox{\rm p.v.}}\limits_{z = a}{ {\Big (} (z-a)^k \cdot
f_{2}(z) } {\Big )} }.
\end{equation}
\end{th}
\begin{co}
Let $f_{1}(z)$ be a holomorphic function at the point $a$ and let $f_{2}(z)$ be a meromorphic function with
simple pole at the same point $a$.$\;$Then
\begin{equation}
\mathop{\mbox{\rm p.v.}}\limits_{z = a}{{\Big (}f_{1}(z) \cdot f_{2}(z){\Big )}}
=
f_{1}(a) \cdot \mathop{\mbox{\rm p.v.}}\limits_{z = a}{f_{2}(z)}
+
{f_{1}\!}^{'}\!(a) \cdot \mathop{\mbox{\rm res}}\limits_{z = a}{f_{2}(z)}.
\end{equation}
The previous formula, in the case of the zeta function $f_2(z) = \zeta(z)$,
is also given~in {\rm \cite{Blau_Visser_Wipf_88_b}}.
\end{co}
The following statement is proved in \cite{Malesevic_03}.
\begin{lm}
\label{TH_Gamma_PV_f}
For the function $f(z)$ with simple pole at point $z = a$ the following is true
\begin{equation}
\label{PV_f_1}
\mathop{\mbox{\rm p.v.}}\limits_{z = a}{f(z)}
=
\lim\limits_{\varepsilon \rightarrow 0}{
\displaystyle\frac{f(a-\varepsilon)+f(a+\varepsilon)}{2}}.
\end{equation}
\end{lm}
\begin{no}
For any meromorphic function $f(z)$, that possesses at worst simple poles,
by formula {\rm (\ref{PV_f_1})} {\em the principal value at point} is defined
as {\em principal part}~{\rm \cite{Blau_Visser_Wipf_88_b}}, which comes
from the quantum field theory {\rm \cite{Blau_Visser_Wipf_88_a}}. Especially,
using formula {\rm (\ref{PV_f_1})}, for the zeta function, {\em the Casimir energy} in physics
is given {\rm \cite{Blau_Visser_Wipf_88_a}}, {\rm \cite{Blau_Visser_Wipf_88_b}}. For a precise
de\-finition of the Casimir energy see {\rm \cite{Elizalde_94}}, {\rm \cite{Elizalde_96}}, {\rm \cite{Elizalde_04}}.
\end{no}
\begin{co}
\label{TH_Gamma_PV}
For gamma function $\Gamma(z)$ it is true {\rm \cite{Slavic_70}},
{\rm \cite{Blau_Visser_Wipf_88_b}}, {\rm \cite{Malesevic_03}}$:$
\begin{equation}
\label{PV_1}
\mathop{\mbox{\rm p.v.}}\limits_{z = -n}{\!\!\Gamma(z)}
=
\lim\limits_{\varepsilon \rightarrow 0}{\!
\displaystyle\frac{\Gamma(-n-\varepsilon)+\Gamma(-n+\varepsilon)}{2}}
=
(-1)^{n}\displaystyle\frac{\Gamma^{'}(n+1)}{\Gamma(n+1)^2}
\quad (n\!\in\!\mbox{\newmat{N}}_{0}).
\end{equation}
\end{co}
\begin{re}
\label{Gamma_Kolicnik_Je_Opadajuci_Niz}
For $n \in \mbox{\newmat{N}}_{0}$ it is true~{\rm \cite{Slavic_70}}$:$

\vspace*{-5.0 mm}

\begin{equation}
\label{H_n}
\displaystyle\frac{\Gamma^{'}(n+1)}{\Gamma(n+1)^2}
=
\displaystyle\frac{\mbox{\small $-\gamma$} + \mbox{\footnotesize $\displaystyle\sum_{k=1}^{n}{\displaystyle\frac{1}{k}}$}}{n!},
\end{equation}

\noindent
where $\gamma$ is Euler's constant.
\end{re}

\noindent
Extension of the domain of the functions $A(z)$ and $A_{1}(z)$, in the sense of
the principal value at point, is given by the following two theorems.
\begin{th}
For alternating Kurepa's function $A(z)$ it is true
\begin{equation}
\label{A_PV}
\mathop{\mbox{\rm p.v.}}\limits_{z = -n}{\!\!A(z)} =\!
\displaystyle\sum\limits_{i=0}^{n-1}{(-1)^{n+1-i}\!\!\!\!\mathop{
\mbox{\rm p.v.}}\limits_{z=-(i-1)}{\!\!\!\!\Gamma(z)}} =\! (-1)^{n+1}
{\bigg (}
1-\!\displaystyle\sum\limits_{i=1}^{n-1}{\displaystyle\frac{\Gamma^{'}(i)}{\Gamma(i)^2}}
{\bigg )} \;\; (n \in \mbox{\newmat{N}}).
\end{equation}
\end{th}
{\bf Proof.} If the equality
$
A(z) = (-1)^{n} A(z+n) + {\big (} \Gamma(z+2) - \ldots + (-1)^{n+1}\Gamma(z+n+1) {\big )}
$
we consider at the point $z = -n$ in the sense of the principal value, on the basis
of (\ref{PV_1}), the equality (\ref{A_PV}) follows. Let us remark that
$\mathop{\mbox{\rm p.v.}}\limits_{z = -1}{\!A(z)} = A(-1) = 1$.~\stop

\medskip
\noindent
The following Ramanujan formula is true:
\begin{equation}
\label{Ramanujan}
\displaystyle\sum\limits_{n=1}^{\infty}{\displaystyle\frac{1 + \frac{1}{2} + \ldots + \frac{1}{n}}{n!} \, x^n}
=
e^x
\displaystyle\sum\limits_{n=1}^{\infty}{\displaystyle\frac{(-1)^{n-1}}{n! \, n} \, x^n},
\end{equation}
for $x \in \mbox{\newmat{C}}$ (see \cite{Berndt_85}, page 46., corollary 2.).
On the basis of the previous formula follows:
\begin{lm}
\label{L_2}
Let us define $L_{2} = \displaystyle\sum\limits_{n=0}^{\infty}{
(-1)^n\!\!\mathop{\mbox{\rm p.v.}}\limits_{z =
-(n-1)}{\!\!\Gamma(z)}}$, then
\begin{equation}
\label{Const_L_2}
L_{2}
=
1 + e \gamma - e \displaystyle\sum\limits_{n=1}^{\infty}{\displaystyle\frac{(-1)^{n-1}}{n! \, n}}
=
1 + e\mbox{\rm Ei}(-1)
\approx
0.403 \, 652 \, 377 \, ,
\end{equation}
where $\mbox{\rm Ei}$ is the function of exponential integral${\,}^{\mbox{\scriptsize $\ast)$}}$\footnote{
$\!\!{}^{\ast)}\,$see formula (\ref{Ei}) in this paper}.
\end{lm}
{\bf Proof.} On the basis of the corollary \ref{TH_Gamma_PV} and the remark \ref{Gamma_Kolicnik_Je_Opadajuci_Niz} we have
\begin{equation}
\label{Const_L_21}
L_{2}
=
1-\displaystyle\sum\limits_{n=0}^{\infty}{\displaystyle\frac{\Gamma^{'}(n\!+\!1)}{\Gamma(n\!+\!1)^2}}
=
1
-
{\bigg (}\displaystyle\sum\limits_{n=1}^{\infty}{\displaystyle\frac{1 + \frac{1}{2} + \ldots + \frac{1}{n}}{n!} }
-
\gamma e{\bigg )}.
\end{equation}
By substitution $x=1$ in formula (\ref{Ramanujan}) we obtain:
\begin{equation}
\label{Const_L_22}
\displaystyle\sum\limits_{n=1}^{\infty}{\displaystyle\frac{1 + \frac{1}{2} + \ldots + \frac{1}{n}}{n!}}
=
e
\displaystyle\sum\limits_{n=1}^{\infty}{\displaystyle\frac{(-1)^{n-1}}{n! \, n}}.
\end{equation}
Next, the following representation of Gompertz
constant $-e \mbox{\rm Ei}(-1)$ (see sequence A073003 in \cite{Sloane_06}) is true

\vspace*{-3.0 mm}

\begin{equation}
\label{Const_L_23}
-
e \mbox{\rm Ei}(-1)
=
e {\bigg (}\displaystyle\sum\limits_{n=1}^{\infty}{\displaystyle\frac{(-1)^{n-1}}{n! \, n}}
-
\gamma {\bigg )}.
\end{equation}
Then, on the basis of (\ref{Const_L_21}), (\ref{Const_L_22}) and (\ref{Const_L_23})
follows (\ref{Const_L_2}).~\stop

\noindent
\begin{th}
\label{A_aux}
For the function $A_{1}(z)$ it is true
\begin{equation}
\mathop{\mbox{\rm p.v.}}\limits_{z = n}{A_{1}(z)}
=
(-1)^n L_{2} + \mathop{\mbox{\rm p.v.}}\limits_{z = n}{A(z)}
\quad (n \!\in\! \mbox{\newmat{Z}}).
\end{equation}
\end{th}
{\bf Proof.} For $n \geq 0$ it is true
\begin{equation}
\mbox{\small $\begin{array}{rcl}
\mathop{\mbox{\rm p.v.}}\limits_{z = n}{A_{1}(z)}
&\!\!=\!\!&
\displaystyle\sum\limits_{i=0}^{\infty}{
(-1)^{i}\!\!\mathop{\mbox{\rm p.v.}}\limits_{z = n+1-i}{\!\Gamma(z)}}
=
(-1)^n
\displaystyle\sum\limits_{i=0}^{\infty}{
(-1)^{i}\!\!\mathop{\mbox{\rm p.v.}}\limits_{z = -(i-1)}{\!\Gamma(z)}}                            \\[1.0 ex]
&\!\!+\!\!&
\displaystyle\sum\limits_{i=1}^{n}{(-1)^{n-i}\Gamma(i+1)}
=
(-1)^{n}L_{2} + A(n).
\end{array}$}
\end{equation}

\noindent
For $n < 0$ it is true
\begin{equation}
\!\!\!\!\!\!\!\mbox{\small $\begin{array}{rcl}
\mathop{\mbox{\rm p.v.}}\limits_{z = n}{A_{1}(z)}
&\!\!\!\!\!=\!\!\!\!\!&
\displaystyle\sum\limits_{i=0}^{\infty}{(-1)^i
\!\!\!\mathop{\mbox{\rm p.v.}}\limits_{z = n+1-i}{\!\!\Gamma(z)}}
= (-1)^{(-n)} {\bigg (}
\displaystyle\sum\limits_{i=(-n)}^{\infty}{(-1)^i
\!\!\!\mathop{\mbox{\rm p.v.}}\limits_{z = -i+1}{\!\!\Gamma(z)}} {\bigg )}                        \\[2.0 ex]
&\!\!\!\!\!=\!\!\!\!\!&
(-1)^{(-n)} {\bigg (}\!
\displaystyle\sum\limits_{i=0}^{\infty}{(-1)^i
\!\!\!\!\mathop{\mbox{\rm p.v.}}\limits_{z = -(i-1)}{\!\!\!\!\Gamma(z)}}\!{\bigg )}
-
(-1)^{(-n)} {\bigg (}\!\displaystyle\sum\limits_{i=0}^{(-n)-1}{\!\!\!(-1)^{-i}
\!\!\!\!\mathop{\mbox{\rm p.v.}}\limits_{z = -(i-1)}{\!\!\!\!\Gamma(z)}}\!{\bigg )}
                                                                               \\[2.0 ex]
&\!\!\!\!\!=\!\!\!\!\!&
(-1)^n L_{2}
+\!
{\bigg (}\!
\displaystyle\sum\limits_{i=0}^{(-n)-1}{
(-1)^{(-n)+1-i} \!\!\!\!\!\mathop{\mbox{\rm p.v.}}\limits_{z = -(i-1)}{\!\!\!\!\Gamma(z)}}
\!{\bigg )}
\!=
(-1)^{n} L_{2} + \mathop{\mbox{\rm p.v.}}\limits_{z = n}{A(z)}.\;\stop
\end{array}$}
\end{equation}

\noindent
\section{\bf \boldmath \hspace*{-7.0 mm}
3 Formula of Slavi\' c's type for alternating Kurepa's function}

\medskip
\noindent
The main result in this section is a new formula of Slavi\' c's type for alternating Kurepa's function
$A(z)$, which is an analog of Slavi\' c's representation of Kurepa's function $K(z)$
\cite{Slavic_73}, \cite{Marichev_83}, \cite{Milovanovic_96}. The following statements are true.
\begin{lm}
\label{A_lema_12}
Function
\begin{equation}
\label{Maple_A_Lema_12}
F(z)
=
\displaystyle\sum\limits_{n=1}^{\infty}{
\!{\bigg (}\displaystyle\sum\limits_{k=2}^{\infty}{
\displaystyle\frac{(-1)^{k-1}(n+k-1)}{(n+k)!} \, z^k {\bigg )}}},
\end{equation}
is entire, whereas the following is true
\begin{equation}
\label{Maple_A_Lema_21}
F(z)
=
\displaystyle\sum\limits_{k=2}^{\infty}{ \!{\bigg (}\displaystyle\sum\limits_{n=1}^{\infty}{
\displaystyle\frac{(-1)^{k-1}(n+k-1)}{(n+k)!} \, z^k {\bigg )}}}
=
-e^{-z}-z+1.
\end{equation}
\end{lm}

\noindent
{\bf Proof.} For $z = 0$ the equality (\ref{Maple_A_Lema_21}) is true. Let us
introduce a sequence of functions
\begin{equation}
\label{Def_f_n}
f_{n}(z)
=
\displaystyle\sum\limits_{k=2}^{\infty}{
\displaystyle\frac{(-1)^{k-1}(n+k-1)}{(n+k)!} \, z^k},
\end{equation}
for $z \in \mbox{\newmat{C}}$ $(n \in \mbox{\newmat{N}})$. Previous series converge over $\mbox{\newmat{C}}$ because,
for $z \neq 0$, it is true that
\begin{equation}
\label{Maple_f_n}
f_{n}(z)
=
\displaystyle\sum\limits_{j=0}^{n+1}{
(-1)^{j+n} {\bigg (}\displaystyle\frac{j}{j!}
\!-\!
\displaystyle\frac{1}{j!} {\bigg )} \, z^{j-n}}
\!+\!
(-1)^{n}e^{-z}(z^{-n+1}\!+\!z^{-n}).
\end{equation}

\vspace*{-2.5 mm}

\noindent
Let us mention that the previous equality is easily checked by the following substi\-tution${}^{C1)}$
\mbox{$e^{-z} = \sum_{k=0}^{\infty}{\frac{(-z)^k}{k!}}$} at the right side of
equality of formula (\ref{Maple_f_n}). Let  $\rho>0$ be fixed. Over the set
\mbox{$\mbox{\newmat{D}}\!=\!\{ z\!\in\!\mbox{\newmat{C}} \, | \, 0\!<\!|z|\!<\!\rho\}$}
let us form an auxiliary function \mbox{$g(z) = (z+1)e^{-z} : \mbox{\newmat{D}} \longrightarrow \mbox{\newmat{C}}$}.
If we denote by $R_{n}(.)$ the remainder of $n$-th order of MacLaurin's expansion, then for $z \in \mbox{\newmat{D}}$
the following representation is true${}^{C2)}$
\begin{equation}
\label{MacLaurin_f_n}
f_{n}(z)
=
\displaystyle\frac{R_{n+1}(g(z))}{(-z)^n}.
\end{equation}
Then, for $|z| < \rho$ it is true
\begin{equation}
\label{Ineq_f_n}
|f_{n}(z)|
\leq
\displaystyle
\frac{e^{\rho}(n+1+\rho)}{(n+2)!}\rho^2.
\end{equation}
Indeed, for $z=0$ the previous inequality is true. Over $\mbox{\newmat{E}} = (0,\rho)$ let us form an auxiliary function
\mbox{$h(t) = (t-1)e^{t} + 2 : \mbox{\newmat{E}} \longrightarrow \mbox{\newmat{R}}^{+}$}. For $z \in \mbox{\newmat{D}}$
and $t = |z| \in \mbox{\newmat{E}}$ there exists $c \in (0,t)$ such that${}^{C3)}$
\begin{equation}
|f_{n}(z)|
\leq
\displaystyle\frac{R_{n+1}(h(t))}{t^n}
=
\displaystyle
\frac{h^{(n+2)}(c)}{(n+2)!}\,t^2
\leq
\displaystyle
\frac{e^{\rho}(n+1+\rho)}{(n+2)!}\rho^2.
\end{equation}
For the function
\begin{equation}
F(z)
=
\displaystyle
\sum\limits_{n=1}^{\infty}{f_{n}(z)},
\end{equation}
it is possible, for $|z|\!<\!\rho$, to apply Weierstrass's double series
Theorem~\cite{Knopp_96} (page~83.). Indeed, on the basis of (\ref{Def_f_n}),
the functions $f_{n}(z)$ are regular for $|z|\!<\!\rho$. On the basis of
(\ref{Ineq_f_n}), the series $\sum_{n=1}^{\infty}{f_{n}(z)}$ is uniformly
convergent for $|z|\!\leq\!r\!<\!\rho$, for every $r\!<\!\rho$. Then on
the basis of the Weierstrass's 's double series Theorem, for $|z|\!<\!\rho$,
the following~is~true
\begin{equation}
F(z)
=
\displaystyle\sum\limits_{k=2}^{\infty}{
\!{\bigg (}\displaystyle\sum\limits_{n=1}^{\infty}{
\displaystyle\frac{(-1)^{k-1}(n+k-1)}{(n+k)!} \, z^k {\bigg )}}}
=
-e^{-z}-z+1,
\end{equation}
because${}^{C4)}$
\begin{equation}
\label{Razmena}
\label{A_sum}
\displaystyle\sum\limits_{n=1}^{\infty}{
\displaystyle\frac{(-1)^{k-1}(n+k-1)}{(n+k)!}} =
\displaystyle\frac{(-1)^{k-1}}{k!}.
\end{equation}
Let us note that $\rho>0$ can be arbitrarily large positive number. Hence, the equality (\ref{Maple_A_Lema_21})
is true for all $z \in \mbox{\newmat{C}}$; i.e. the function $F(z)$ is entire. \stop
\begin{lm}
\label{A_lema_3}
For $z \in \mbox{\newmat{C}}$ it is true
\begin{equation}
\label{A_sum_sum} (z\!+\!1)
\displaystyle\sum\limits_{n=1}^{\infty}{
\displaystyle\sum\limits_{k=1}^{\infty}{
\displaystyle\frac{(-1)^{k}}{(k\!+\!n)!} \, z^{k}}}
=
-e^{-z}+(1\!-\!e)\,z+1.
\end{equation}
\end{lm}

\noindent
{\bf Proof.} On the basis of the Lemma \ref{A_lema_12} it is true that
\begin{equation}
\!\!\mbox{\small $
\begin{array}{rcl}
(z\!+\!1)
\displaystyle\sum\limits_{n=1}^{\infty}{
\displaystyle\sum\limits_{k=1}^{\infty}{
\displaystyle\frac{(-z)^{k}}{(k\!+\!n)!}}}
&\!\!=\!\!&
\displaystyle\sum\limits_{n=1}^{\infty}{
\displaystyle\sum\limits_{k=1}^{\infty}{
\displaystyle\frac{(-1)^{k}}{(k\!+\!n)!}\,(z^{k+1}\!+z^{k})}}                  \\[2.5 ex]
&\!\!=\!\!&
\displaystyle\sum\limits_{n=1}^{\infty}{
\!{\bigg (}\!\!-\displaystyle\frac{z}{(n\!+\!1)!}
+\!\displaystyle\sum\limits_{k=2}^{\infty}{
{\Big (}\displaystyle\frac{(-1)^{k-1}}{(k\!+\!n\!-\!1)!}
-\displaystyle\frac{(-1)^{k-1}}{(k\!+\!n)!}{\Big )}z^{k}}{\bigg )}}            \\[2.5 ex]
&\!\!\mathop{=}\limits_{(\ref{Maple_A_Lema_21})}\!\!&
\displaystyle\sum\limits_{k=2}^{\infty}{
{\bigg (}\displaystyle\sum\limits_{n=1}^{\infty}{
\displaystyle\frac{(-1)^{k-1}(n\!+\!k\!-\!1)}{(n\!+\!k)!}}{\bigg )}\,z^{k}}
-\displaystyle\sum\limits_{n=1}^{\infty}{
\displaystyle\frac{z}{(n\!+\!1)!}}                                            \\[5.0 ex]
&\!\!\mathop{=}\limits_{(\ref{Maple_A_Lema_21})}\!\!&
-e^{-z}+(1\!-\!e)\,z+1.\;\stop
\end{array}$}
\end{equation}
\begin{th}
For alternating Kurepa's function $A(z)$ the following representation is true
\begin{equation}
\label{A_sum_Slavic}
A(z)
=
-
{\big (}1 + e \, \mbox{\rm Ei}(-1){\big )}(-1)^z
+
\displaystyle\frac{\pi e}{\sin \pi z}
+
\displaystyle\sum\limits_{n=0}^{\infty}{(-1)^{n}\Gamma(z+1-n)},
\end{equation}
where the values in the previous formula, in integer points $z$,
are determined in the sense of the principal value at point.
\end{th}

\noindent
{\bf Proof.} For $-(n\!+\!1)\!<\!\mbox{\rm Re}\,z\!<\!-n$ and $n=0,1,2,\ldots$ the following
formula is true~\cite{Bateman_65}:
\begin{equation}
\Gamma(z)
=\!\!
\displaystyle\int\limits_{0}^{+\infty}\!{{\bigg (}\!
e^{-t} - \displaystyle\sum\limits_{m=0}^{n}{
\displaystyle\frac{(-t)^{m}}{m!}}\!{\bigg )}t^{z-1}\:dt}.
\end{equation}
Hence, for $0 < \mbox{\rm Re}\:z < 1$ and $n=1,2,\ldots$ the following formula is true
\begin{equation}
\Gamma(z-n)
=\!\!
\displaystyle\int\limits_{0}^{+\infty}\!{{\bigg (}\!
e^{-t} - \displaystyle\sum\limits_{m=0}^{n-1}{
\displaystyle\frac{(-t)^{m}}{m!}}\!{\bigg )}t^{z-n-1}\:dt}.
\end{equation}
Further we observe the following difference
\begin{equation}
\begin{array}{rcl}
A(z) - \displaystyle\sum\limits_{n=0}^{+\infty}{(-1)^{n}\Gamma(z\!+\!1\!-\!n)}
\!\!&\!\!=\!\!&\!\!
\!\displaystyle\int\limits_{0}^{+\infty}{\!
e^{-t}\displaystyle\frac{t^{z+1}\!-\!(-1)^z t}{t\!+\!1}\:dt}
\;-\!
\displaystyle\int\limits_{0}^{+\infty}{\!
e^{-t}t^{z}\:dt}                                                               \\[2.0 ex]
\!\!&\!\!-\!\!&
\!\displaystyle\sum\limits_{n=1}^{\infty}{
(-1)^{n}\!\!\displaystyle\int\limits_{0}^{+\infty}{\!\!{\bigg (}\!
e^{-t}\!-\!\displaystyle\sum\limits_{m=0}^{n}{\!
\displaystyle\frac{(-t)^{m}}{m!}}\!{\bigg )}t^{z\!-\!n}\:dt}}.
\end{array}
\end{equation}
For $0 < \mbox{\rm Re}\:z < 1$ the following derivation is true
\begin{equation}
\label{A_derive}
\mbox{\small $\begin{array}{rl}
& A(z) - \displaystyle\sum\limits_{n=0}^{+\infty}{(-1)^{n}\Gamma(z\!+\!1\!-\!n)}
                                                                               \\[2.0 ex]
&
= -\!\displaystyle\int\limits_{0}^{\infty}{\!e^{-t}
\displaystyle\frac{t^{z}\!-\!(-1)^{z-1}t}{t\!+\!1} \: dt}
-\!
\displaystyle\int\limits_{0}^{+\infty}{\!
\displaystyle\sum\limits_{n=1}^{\infty}{\!
(-1)^{n}\!{\bigg (}\!e^{-t}\!-\!\displaystyle\sum\limits_{m=0}^{n}{\!
\displaystyle\frac{(-t)^{m}}{m!}}\!{\bigg )}t^{z\!-\!n}\:dt}}                  \\[2.0 ex]
& =
\displaystyle\int\limits_{0}^{\infty}{\!{\bigg (}
\!(-1)^{z-1}e^{-t}\displaystyle\frac{t}{t\!+\!1}
-
e^{-t}\displaystyle\frac{t^{z}}{t\!+\!1}
-
\displaystyle\sum\limits_{n=1}^{\infty}{\!
(-1)^n \!\!\! \displaystyle\sum\limits_{m=n+1}^{\infty}{
\displaystyle\frac{(-t)^{m}}{m!}}\,t^{z\!-\!n}\!{\bigg )}\,dt}}                \\[2.0 ex]
& =
\displaystyle\int\limits_{0}^{\infty}{\!{\bigg (}
\!(-1)^{z-1}e^{-t}\displaystyle\frac{t}{t\!+\!1}
-
\displaystyle\frac{t^{z}}{t\!+\!1}
\!{\bigg (}\!e^{-t}+(t\!+\!1)\!
\mathop{\displaystyle\sum}\limits_{n\!\,=\!\,1}^{\infty}{
\mathop{\!\!\displaystyle\sum}\limits_{\,\,m\!\,=\!\,n+1}^{\infty}{\!\!\!
\displaystyle\frac{(-1)^{m+n}}{m!}\,t^{m-n}}}\!{\bigg )}\!\!{\bigg )}\,dt}     \\[2.0 ex]
& =
\displaystyle\int\limits_{0}^{\infty}{\!{\bigg (}
(-1)^{z-1}e^{-t}\displaystyle\frac{t}{t\!+\!1}
-
\displaystyle\frac{t^{z}}{t\!+\!1}
\!{\bigg (}\!e^{-t}+(t\!+\!1)\!
\displaystyle\sum\limits_{n=1}^{\infty}{
\displaystyle\sum\limits_{k=1}^{\infty}{
\displaystyle\frac{(-t)^{k}}{(k\!+\!n)!}}}\!{\bigg )}\!\!{\bigg )}dt}          \\[2.0 ex]
& \mathop{=}\limits_{(\ref{A_sum_sum})}
\displaystyle\int\limits_{0}^{\infty}{\!{\bigg (}
(-1)^{z-1}e^{-t}\displaystyle\frac{t}{t\!+\!1}
+
\displaystyle\frac{t^{z}}{t+1} {\Big (} (e-1)\,t - 1 {\Big )}\!\!{\bigg )} \, dt}.
\end{array}$}
\end{equation}
Integral at the right side of the equality (\ref{A_derive}), which converges in ordinary sense,
will be substituted by the sum of two integrals which converge in the ordinary sense too.
Namely, using the function of exponential integral, i.e. formula 8.211-1 \cite{Rizik_71}:
\begin{equation}
\label{Ei}
\mbox{\rm Ei}(x) =
\displaystyle\int\limits_{-\infty}^{x}{\!\displaystyle\frac{e^t}{t} \, dt}
\;\; (x<0),
\end{equation}
and using the formulas 3.351-5 and 3.241-2 from \cite{Rizik_71}:
\begin{equation}
\displaystyle\int\limits_{0}^{\infty}{
\displaystyle\frac{e^{-t}}{t+1} \, dt}
=
- e \, \mbox{Ei}(-1)
\quad \; \mbox{and} \; \quad
\displaystyle\int\limits_{0}^{\infty}{
\displaystyle\frac{t^{z-1}}{t+1} \, dt}
=
\displaystyle\frac{\pi}{\sin \pi z}
\end{equation}
\noindent
we can conclude that formula (\ref{A_sum_Slavic}) is true for $0 < \mbox{\rm Re}\:z < 1$.
According to Riemann's Theorem we can conclude that formula (\ref{A_sum_Slavic})
is true for each complex $z$. Namely, formula (\ref{A_sum_Slavic}), in integer points $z$,
is true in the sense of Cauchy's principal value at point on the basis of the Lemma \ref{L_2}
and the Theorem \ref{A_aux}.$\;\stop$

\begin{co}
For alternating Kurepa's function $A(z)$ the following representation is true
\begin{equation}
\label{A_sum_Slavic_2}
\quad
A(z)
=
{\Big (} e \displaystyle\sum\limits_{n=1}^{\infty}{\displaystyle\frac{(-1)^{n-1}}{n! \, n}} - 1 - e \gamma {\Big )}(-1)^z
+
\displaystyle\frac{\pi e}{\sin \pi z}
+
\displaystyle\sum\limits_{n=0}^{\infty}{(-1)^{n}\Gamma(z+1-n)},
\end{equation}
where the values in the previous formula, in integer points $z$,
are determined in the sense of the principal value at point.
\end{co}
\begin{co}
Function $A_{1}(z)$ is a meromorphic function with simple poles in
integer points $z=m$ $(m\!\in\!\mbox{\newmat{Z}})$ and with residue values
\begin{equation}
\mathop{\mbox{\rm res}}\limits_{z = m}{A_{1}(z)}
=
(-1)^{m-1}e + \mathop{\mbox{\rm res}}\limits_{z = m}{A(z)} \quad (m\!\in\!\mbox{\newmat{Z}}).
\end{equation}
At the point $z = \infty$ function $A_{1}(z)$ has an essential singularity.
\end{co}

\noindent
\section{\bf \boldmath \hspace*{-7.0 mm}
4 Some representations of functions $A(z)$ and $A_{1}(z)$ via incomplete gamma function}

In this section we give some representations of functions $A(z)$ and $A_{1}(z)$ via
gamma and incomplete gamma functions, where the last ones 
are defined by integrals:
\begin{equation}
\gamma(a,z)
=
\displaystyle\int\limits_{0}^{z}{e^{-t} t^{\alpha-1} \, dt}
\quad \; \mbox{and} \; \quad
\Gamma(a,z)
=
\displaystyle\int\limits_{z}^{\infty}{e^{-t} t^{\alpha-1} \, dt}.
\end{equation}
Parameters $\alpha$ and $z$ are complex numbers and $t^\alpha$ takes its principal value.
Let us remark that the value $\gamma(\alpha,z)$ exists for $\mbox{\rm Re} \, \alpha > 0$
and the value $\Gamma(\alpha,z)$ exists for $|\mbox{\rm arg} \, z| < \pi$. Then,  we have:
$\gamma(a,z) + \Gamma(a,z) = \Gamma(a)$. Analytical continuation can be obtained on
the basis of representation of the $\gamma$ function using series.

\medskip
\noindent
On the basis of the well-known formula 13., page 325., from \cite{Prudnikov_81}:
\begin{equation}
\label{Temme_1}
\displaystyle\int\limits_{0}^{\infty}{e^{-t} \displaystyle\frac{t^{z+1}}{t+1} \, dt}
=
e \Gamma(z+2)\Gamma(-z-1,1)
\end{equation}
we directly get some representations of functions $A(z)$ and $A_{1}(z)$ via
an incomplete gamma function.
\begin{th}
For functions $A(z)$ and $A_{1}(z)$ the following representations are true
\begin{equation}
\label{Repr_A}
A(z)
=
-{\big (}1\!+\!e \, \mbox{\rm Ei}(-1){\big )} (-1)^z
+
e \, \Gamma(z+2) \, \Gamma(-\!z\!-\!1,1)
\end{equation}
and

\medskip

\noindent
\begin{equation}
\label{Repr_A1}
A_{1}(z)
=
-
\displaystyle\frac{\pi e}{\sin \pi z}
+
e \, \Gamma(z+2) \, \Gamma(-\!z\!-\!1,1),
\end{equation}

\smallskip

\noindent
where the values in the previous formula, in integer points $z$,
are determined in the sense of the principal value at point.
\end{th}
\begin{re}
Formula {\rm (\ref{Repr_A})} also is given in {\rm \cite{Petojevic_02}}.
\end{re}

\break

\noindent
\section{\bf \boldmath \hspace*{-7.0 mm}
5 Differential transcendency of functions $A(z)$ and $A_{1}(z)$}

In this section we provide one statement about differential transcendency of some
solutions of functional equation (\ref{A_FE_1}). Namely, using the method for
proving~of~the differential transcendency from papers
\cite{MijajlovicMalesevic_06} and \cite{MijajlovicMalesevic_07}
we can conclude that the following statement is true:
\begin{th}
Let $\mbox{$\cal M$}_{\mbox{\newmatsm{D}}}$ be a differential field of the meromorphic functions over
a domain $\mbox{\newmat{D}} \!\subseteq\! \mbox{\newmat{C}} \backslash \mbox{\newmat{Z}}^{-}$.
If $g \!=\! g(z) \!\in\! {\cal M}_{\mbox{\newmatsm{D}}}$ is one solution of the functional \mbox{equation}~\mbox{\rm
(\ref{A_FE_1})}, then $g$ is not a solution of any algebraic-differential equation over the field of rational functions
$\mbox{\newmat{C}}(z)$.
\end{th}
\begin{co}
\label{Posledica_A_K1} Especially the functions $A(z)$ and $A_{1}(z)$ are not solutions of
any algebraic-differential equation over the field of rational functions $\mbox{\newmat{C}}(z)$.
\end{co}

\bigskip

\rightline{Received$\,:\;8/2/2008$}

\newpage

$\,$

\bigskip

{\bf Some additional comments (outside the official version)}

\bigskip

\bigskip

\bigskip

{ 

\footnotesize

\noindent Equation (25)

\noindent
---------------------------------------------

\noindent ${}^{C1)}$

\vspace*{-5.0 mm}

$$
\begin{array}{l}
\;\;\;\displaystyle\sum\limits_{j=0}^{n+1}{ (-1)^{j+n} {\bigg
(}\displaystyle\frac{j}{j!} \!-\! \displaystyle\frac{1}{j!} {\bigg
)} \, z^{j-n}} \!+\!
(-1)^{n}e^{-z}(z^{-n+1}\!+\!z^{-n})                                             \\[1.0 ex]
= \displaystyle\sum\limits_{j=0}^{n+1}{ (-1)^{j+n} {\bigg
(}\displaystyle\frac{j}{j!} \!-\! \displaystyle\frac{1}{j!} {\bigg
)} \, z^{j-n}} \!+\! (-1)^{n} {\bigg
(}\displaystyle\sum\limits_{j=0}^{\infty}{\displaystyle\frac{(-1)^{j}}{j!}z^{j}}{\bigg
)}
(z^{-n+1}\!+\!z^{-n})                                                          \\[1.0 ex]
= \displaystyle\sum\limits_{j=0}^{n+1}{(-1)^{j+n}
\displaystyle\frac{j}{j!} \, z^{j-n}} \!-\!
\displaystyle\sum\limits_{j=0}^{n+1}{(-1)^{j+n}
\displaystyle\frac{1}{j!} \, z^{j-n}} \!+\!
(-1)^{n}\!\!\displaystyle\sum\limits_{j=0}^{\infty}{\displaystyle\frac{(-1)^{j}}{j!}z^{j+1-n}}
\!+\!
(-1)^{n}\!\!\displaystyle\sum\limits_{j=0}^{\infty}{\displaystyle\frac{(-1)^{j}}{j!}z^{j-n}}
                                                                                \\[1.0  ex]
=
(-1)^{n}\!\!\!\!\displaystyle\sum\limits_{j=n+1}^{\infty}{\displaystyle\frac{(-1)^{j}}{j!}z^{j+1-n}}
+
(-1)^{n}\!\!\!\!\displaystyle\sum\limits_{j=n+2}^{\infty}{\displaystyle\frac{(-1)^{j}}{j!}z^{j-n}}
                                                                                \\[1.0  ex]
=
(-1)^{n}\!\displaystyle\sum\limits_{k=2}^{\infty}{\displaystyle\frac{(-1)^{n+k-1}}{(n+k-1)!}z^{k}}
+
(-1)^{n}\!\displaystyle\sum\limits_{k=2}^{\infty}{\displaystyle\frac{(-1)^{n+k}}{(n+k)!}z^{k}}
                                                                                \\[1.0  ex]
= \displaystyle\sum\limits_{k=2}^{\infty}{(-1)^{k-1}{\Big
(}\displaystyle\frac{1}{(n+k-1)!} -
\displaystyle\frac{1}{(n+k)!}{\Big )} z^{k}}                                    \\[1.0 ex]
=
\displaystyle\sum\limits_{k=2}^{\infty}{\displaystyle\frac{(-1)^{k-1}(n+k-1)}{(n+k)!}
\, z^k} \quad (z \neq 0) \; .
\end{array}
$$}

\bigskip

\bigskip

\bigskip

{ 

\footnotesize

\noindent Equation (26)

\noindent
---------------------------------------------

\noindent ${}^{C2)}$

\vspace*{-5.0 mm}

$$
\begin{array}{rcl}
g(z) \!\!&\!\!=\!\!&\!\! (z+1) \, e^{-z} = (z+1)
\displaystyle\sum\limits_{j=0}^{\infty}{\displaystyle\frac{(-1)^{j}}{j!}z^{j}}
=
\displaystyle\sum\limits_{j=0}^{\infty}{\displaystyle\frac{(-1)^{j}}{j!}z^{j+1}}
+
\displaystyle\sum\limits_{j=0}^{\infty}{\displaystyle\frac{(-1)^{j}}{j!}z^{j}}  \\[1.0 ex]
\!\!&\!\!=\!\!&\!\!
\displaystyle\sum\limits_{j=1}^{\infty}{\displaystyle\frac{(-1)^{j-1}}{(j\!-\!1)!}z^{j}}
+
\displaystyle\sum\limits_{j=0}^{\infty}{\displaystyle\frac{(-1)^{j}}{j!}z^{j}}
= 1 + \displaystyle\sum\limits_{j=1}^{\infty}{ (-1)^{j-1} {\Big (}
\displaystyle\frac{1}{(j\!-\!1)!} - \displaystyle\frac{1}{j!} {\Big
)}
z^{j}}                                                                          \\[1.0 ex]
\!\!&\!\!=\!\!&\!\! 1 +
\displaystyle\sum\limits_{j=1}^{\infty}{(-1)^{j-1}
\displaystyle\frac{j-1}{j!} z^{j}} {\Big (} \! = 1 \!-\!
\displaystyle\frac{1}{2}z^2 \!+\! \displaystyle\frac{1}{3}z^3 \!-\!
\displaystyle\frac{1}{8}z^4 \!+\! \displaystyle\frac{1}{30}z^5 \!-\!
\displaystyle\frac{1}{144}z^6 \!+\! \ldots {\Big )},
\end{array}
\leqno (p)
$$
$$
R_{n+1}{\big (} g(z) {\big )} \;\; \mathop{=}\limits_{(p)} \;
\displaystyle\sum\limits_{j=n+2}^{\infty}{
(-1)^{j-1}\displaystyle\frac{j-1}{j!}z^{j}} \;\; = \;\;
\displaystyle\sum\limits_{k=2}^{\infty}{
(-1)^{n+k-1}\displaystyle\frac{n+k-1}{(n+k)!}z^{n+k}} \;\; ; \leqno
(q)
$$
then:
$$
f_{n}(z) \;\; = \;\; \displaystyle\sum\limits_{k=2}^{\infty}{
(-1)^{k-1}\displaystyle\frac{n+k-1}{(n+k)!}z^{k}} \;\;
\mathop{=}\limits_{(q)} \; \displaystyle\frac{R_{n+1}{\big
(}g(z){\big )}}{(-z)^{n}} \quad (z \!\in\! \mbox{\newmatmb{D}}) \; .
$$}

\break


{ 

\footnotesize

\noindent Inequalities (27), (28)

\noindent
---------------------------------------------

\noindent ${}^{C3)}$

\vspace*{-5.0 mm}

$$
\begin{array}{rcl}
h(t) \!\!&\!\!=\!\!&\!\! (t\!-\!1) \, e^{t} + 2 \; = \; (t\!-\!1) \!
\displaystyle\sum\limits_{j=0}^{\infty}{\displaystyle\frac{1}{j!}t^{j}}
+ 2 \; =
\displaystyle\sum\limits_{j=0}^{\infty}{\displaystyle\frac{1}{j!}t^{j+1}}
\!-\!
\displaystyle\sum\limits_{j=0}^{\infty}{\displaystyle\frac{1}{j!}t^{j}} + 2     \\[1.0 ex]
\!\!&\!\!=\!\!&\!\! 2 +
\displaystyle\sum\limits_{j=1}^{\infty}{\displaystyle\frac{1}{(j\!-\!1)!}t^{j}}
\;-\;
\displaystyle\sum\limits_{j=0}^{\infty}{\displaystyle\frac{1}{j!}t^{j}}
\;=\; 1 + \displaystyle\sum\limits_{j=1}^{\infty}{ {\Big (}\,
\displaystyle\frac{1}{(j-1)!} \,-\, \displaystyle\frac{1}{j!}
\;{\Big )}
t^{j}}                                                                          \\[1.0 ex]
\!\!&\!\!=\!\!&\!\! 1 +
\displaystyle\sum\limits_{j=1}^{\infty}{\displaystyle\frac{j-1}{j!}
t^{j}} {\Big (} \!\! = 1 \!+\! \displaystyle\frac{1}{2}t^2 \!+\!
\displaystyle\frac{1}{3}t^3 \!+\! \displaystyle\frac{1}{8}t^4 \!+\!
\displaystyle\frac{1}{30}t^5 \!+\! \displaystyle\frac{1}{144}t^6
\!+\! \ldots {\Big )} \, ,
\end{array}
\leqno (\ppsm{p})
$$
$$
R_{n+1}{\big (} h(t) {\big )} \;\mathop{=}\limits_{(\ppssm{p})}
\displaystyle\sum\limits_{j=n+2}^{\infty}{
\displaystyle\frac{j-1}{j!}t^{j}} =
\displaystyle\sum\limits_{k=2}^{\infty}{
\displaystyle\frac{n+k-1}{(n+k)!}t^{n+k}} \; , \leqno (\ppsm{q})
$$
$$
h^{(k)}(t) = e^{t} (t + k - 1) \quad (k=2,3,\ldots) \, ; \leqno
(\ppsm{r})
$$

\medskip
\noindent then, for $z \!\in\! \mbox{\newmatmb{D}}$ and $t \!=\! |z| \!\in\! \mbox{\newmatmb{E}} \!=\! (0,\rho)$
the following inequalities are true:
$$
\begin{array}{rcl}
|f_{n}(z)| \!\!&\!\! \leq \!\!&\!\!
\displaystyle\sum\limits_{k=2}^{\infty}{\displaystyle\frac{n+k-1}{(n+k)!}
\, |z|^{k}} \;=\;
\displaystyle\sum\limits_{k=2}^{\infty}{\displaystyle\frac{n+k-1}{(n+k)!}
\, t^{k}} \;\;\mathop{=}\limits_{(\ppssm{q})}\;\;
\displaystyle\frac{R_{n+1}{\big (}h(t){\big )}}{t^n}                            \\[1.0 ex]
\!\!&\!\!=\!\!&\!\! \displaystyle\frac{{\big
(}h^{(n+2)}(c)/(n+2)!{\big )} \, t^{n+2}}{t^n}
\;\; {\big (}\mbox{by the Lagrange remainder theorem, for some}\;c\!\in\!(0,t) {\big )}      \\[2.0 ex]
\!\!&\!\!=\!\!&\!\! \displaystyle\frac{h^{(n+2)}(c)}{(n+2)!} \, t^2
\;\mathop{=}\limits_{(\ppssm{r})}\;
\displaystyle\frac{e^c(c+n+1)}{(n+2)!} \, t^2 \;\leq\;
\displaystyle\frac{e^{\rho}(\rho + n + 1)}{(n + 2)!} \, \rho^2 \, .
\end{array}
$$}

\bigskip

\bigskip

\bigskip

{ 

\footnotesize

\noindent Equation (31)

\noindent
---------------------------------------------

\noindent ${}^{C4)}$

\vspace*{-5.0 mm}
$$
\begin{array}{rcl}
\displaystyle\sum\limits_{n=1}^{\infty}{\displaystyle\frac{(-1)^{k-1}(n+k-1)}{(n+k)!}}
\!\!&\!\!=\!\!&\!\!
(-1)^{k-1}\!\displaystyle\sum\limits_{n=1}^{\infty}{ {\Big (} \,
\displaystyle\frac{n+k}{(n+k)!}-\displaystyle\frac{1}{(n+k)!}
\,{\Big )}}                                                                     \\[1.0 ex]
\!\!&\!\!=\!\!&\!\!
(-1)^{k-1}\!\displaystyle\sum\limits_{n=1}^{\infty}{ {\Big (}
\displaystyle\frac{1}{(n\!+\!k\!-\!1)!}-\displaystyle\frac{1}{(n\!+\!k)!}
{\Big )}} \;=\; \displaystyle\frac{(-1)^{k-1}}{k!}.
\end{array}
$$}

\end{document}